\documentclass[12pt]{article}
\usepackage{amssymb,amsmath,amsthm,tikz-cd,amscd,cite,hyperref, pgfplots, authblk, tikz, caption, graphicx}
\usepackage[letterpaper, margin=1in]{geometry}

\title{Heat Kernel Estimates for Schr\"odinger Operators in the Domain Above a Bounded Lipschitz Function}

\author[1]{Anthony Graves-McCleary}
\affil[1]{Department of Mathematics, Cornell University, Ithaca, NY, USA ag2537@cornell.edu}
\date{}
\newcommand{\set}[1]{\left\{#1\right\}}
\newcommand{\Real}{\mathbf R}

\newcommand{\N}{\mathbf{N}}

\newcommand{\C}{\mathbf{C}}

\newcommand{\brac}[1]{\left< #1 \right>}
\newcommand{\abs}[1]{\left\vert#1\right\vert}

\newcommand{\fct}[3]{#1 \colon #2 \rightarrow #3}

\renewcommand{\phi}{\varphi}

\newcommand{\mc}{\mathcal}

\newcommand{\grad}{\nabla}

\newtheorem{theorem}{Theorem}[section]
\newtheorem{proposition}[theorem]{Proposition}
\newtheorem{lemma}[theorem]{Lemma}
\newtheorem{definition}[theorem]{Definition}
\newtheorem{corollary}[theorem]{Corollary}

\usetikzlibrary{arrows, matrix}

\numberwithin{equation}{section}

\begin{document}
\maketitle

\begin{abstract}
We give matching upper and lower bounds for the Dirichlet heat kernel of a Schr\"odinger operator $\Delta+W$ in the domain above the graph of a bounded Lipschitz function, in the case when $W$ decays away from the boundary faster than quadratically.
\end{abstract}

\section{Introduction and Overview}

We are interested in estimating the heat kernel of operators of the form $\Delta+W$ in a domain $\Omega$ in $\Real^N$, where $\Delta=-\sum_{j=1}^N \frac{\partial^2}{\partial x_j^2}$ is the Laplacian and we impose Dirichlet boundary conditions. 

The specific setting we study is as follows. Let $\fct{f}{\Real^{N-1}}{\Real}$ be bounded and Lipschitz, with $N\geq 2$. Let $$\Omega=\set{(x_1, \dots, x_{N-1}, x_N)\colon x_N>f(x_1, \dots, x_{N-1})}$$ be the domain above the graph of $f$. Let $\fct{W}{\Omega}{\Real}$ be smooth such that there exist $c, \epsilon>0$ with $$0\leq W(x)\leq c(1+\delta(x))^{-(2+\epsilon)}$$ for all $x\in \Omega$, where $\delta(x)=\inf\set{\abs{x-y}\colon y\in \partial \Omega}$ is the distance to the boundary. We are interested in the following problem:\\

Find matching upper and lower bounds for the Dirichlet heat kernel of $\Delta+W$ in $\Omega$.\\ 

We give the answer to this problem in our main result, Theorem \ref{maithm}.\\

In this section we give an overview of related work. Song \cite{song1} gave bounds on the heat kernel of the Laplacian in the domain above the graph of a bounded $C^{1,1}$ function. Hirata in \cite{hirata1} gave estimates on the profile of certain Schr\"odinger operators in uniform cones, on which we make two observations. First, that some uniform cones can be realized as domains above the graph of a Lipschitz function, although this function is usually unbounded. Second, the Schr\"odinger operators studied by Hirata are of the form $(1+\abs{x})^{p-q}\delta(x)^{-p}$ for $p<2<q$ where $\delta$ is distance from the boundary of the domain. Notably, these may be unbounded near the boundary and they decay at infinity. The potentials studied in our paper are bounded, but decay only with distance from the boundary, not the origin.\\

Chen, Kim, and Song in \cite{chenkimsong1} gave estimates for the Green's function of some non-local operators in half-space-like $C^{1,1}$ open sets, and later Chen and Song in \cite{chensong1} gave heat kernel estimates in a similar setting. Gyrya and Saloff-Coste in \cite{gyryasaloffcoste1} gave heat kernel estimates for the Laplacian in the domain above a Lipschitz graph. We will discuss this result in detail later, as it is of particular importance to our results.

\subsection{Notation}

Given $x\in \Real^N$ we let $\abs{x}=\left(\sum_{j=1}^Nx_j^2\right)^{1/2}$ denote the Euclidean norm. We let $B(x,r)=\set{y \in \Real^N\colon \abs{x-y}<r}$ denote the open ball of radius $r>0$ and center $x$ in $\Real^N$. The Laplacian is given by $\Delta=-\sum_{j=1}^N \frac{\partial^2}{\partial x_j^2}$; we caution that this sign convention is the opposite of some authors.\\

Given a domain $\Omega\subseteq \Real^N$, we let $C^\infty(\Omega)$ denote the set of smooth functions on $M$. For $x\in \Omega$, we let $\delta(x)=d(x,\partial\Omega)=\inf\set{\abs{x-y}\colon y\notin \Omega}$ denote distance to the boundary. We define $\brac{x}:=1+\delta(x)$.\\

 Given quantities $a$ and $b$, we let $a\preceq b$ denote the existence of a constant $C>0$ such that $a\leq Cb$. We let $a\asymp b$ denote that $a\preceq b$ and $b\preceq a$. Exact values of constants may sometimes change from line to line.\\

Given $x=(x_1, \dots, x_N)\in \Real^N$, we let $x'=(x_1, \dots, x_{N-1})$ denote the projection of $x$ in $\Real^{N-1}$. In practice we identify the ordered pair $(x', x_N)$ with $x$. Given $r>0$, we define $x+r:=(x', x_N+r)$. The expressions $dx$, $dy$, $dy'$, $dy_N$, etc. will denote Lebesgue measure of the appropriate dimension.\\

\section{Uniform Domains and the Boundary Harnack Principle}

In this section we discuss uniform domains in Euclidean space, with an emphasis on the Boundary Harnack Principle.

\begin{definition}
Let $\Omega\subseteq \Real^N$ be a domain. We say that $\Omega$ is a \textbf{uniform domain} if there exist $c_0, C_0>0$ such that for any $x,y\in \Omega$, there exists a rectifiable curve $\gamma$ of length at most $C_0$ joining $x$ to $y$ and satisfying, for all $z\in \gamma([0,1])$, $$d(x,\partial \Omega)\geq c_0\frac{\abs{x-z}\abs{y-z}}{\abs{x-y}}.$$
\end{definition}

\begin{definition}
Given a domain $\Omega\subseteq \Real^N$, we define the intrinsic distance between $x,y\in \Omega$, denoted $d(x,y)$, to be the infimum of the lengths of rectifiable curves joining $x$ to $y$.
\end{definition}

\noindent\textbf{Remark:} If $\Omega$ is a uniform domain, then it is straightforward to verify that the instrinsic distance $d(x,y)$ between $x,y\in \Omega$ is comparable to the Euclidean distance $\abs{x-y}$. We will use this later.\\

The domain above the graph of a Lipschitz function is a uniform domain, as seen in e.g. Gyrya and Saloff-Coste \cite{gyryasaloffcoste1} . For us, the most important property of a general uniform domain is the Boundary Harnack Principle, which we introduce below.\\

\begin{theorem}\label{boundaryharnack}
Let $\Omega$ be an unbounded uniform domain in $\Real^N$. There exist constants $A_0, A_1>1$ such that for any $\xi\in \partial \Omega$, any $r>0$, and any positive harmonic functions $u$ and $v$ in $B(\xi, A_0r)\cap \Omega$ with Dirichlet boundary conditions along $\partial \Omega$, we have that for all $x,y\in B(\xi, r)\cap \Omega$, $$\frac{u(x)}{u(y)}\leq A_1\frac{v(x)}{v(y)}.$$
\end{theorem}

\noindent\textbf{Proof:} See e.g. Gyrya and Saloff-Coste \cite{gyryasaloffcoste1}, Theorem 4.2 for an even more general version.\qed\\

\noindent\textbf{Remark:} The Boundary Harnack Principle has a rich history. It originated in papers of Ancona \cite{ancona1}, Dahlberg \cite{dahlberg1}, and Wu \cite{wu1} in the late 1970s. See the textbook by Armitage and Gardiner \cite{armitagegardiner1} for a presentation in Lipschitz domains.\\

\begin{definition}
Let $\Omega$ be an unbounded uniform domain in $\Real^N$ and let $W\in C^\infty(\Omega)$. We say that a function $h$ on $\Omega$ is a \textbf{profile} of $\Delta+W$ on $\Omega$ if $h>0$ on $\Omega$, $h$ satisfies Dirichlet boundary conditions along $\partial \Omega$ and $$(\Delta+W)h=0.$$ If $W=0$ we say that $h$ is a \textbf{harmonic profile} on $\Omega$.
\end{definition}

\begin{lemma}
Let $\Omega$ be an unbounded uniform domain $\Real^N$. Then there exists a harmonic profile $h>0$ on $\Omega$. Furthermore, if $g,h>0$ are harmonic profiles on $\Omega$, then there exists $c>0$ such that $g=ch$.
\end{lemma}

\noindent\textbf{Proof:} The proof of existence can be found in e.g. Gyrya and Saloff-Coste \cite{gyryasaloffcoste1}, Theorem 4.16. We will discuss some of the details of this proof later in our paper, as they will be useful in the proof of Theorem \ref{profilebound}.\\

For uniqueness up to a constant, a much more general version of this argument is given in e.g. Kajino and Murugan \cite{kajinomurugan1}, Lemma 3.19.\qed\\

We finish this section with a useful result about the domain above the graph of a Lipschitz function.

\begin{proposition}\label{volumeest}
Let $\Omega\subseteq \Real^N$ be the domain above the graph of a Lipschitz function, and let $h>0$ be a harmonic profile on $\Omega$. Then $h$ is non-decreasing in the vertical direction, and if $x\in \Omega$ and $r>0$, if $x=(x', x_N)$ and $x+r:=(x', x_N+r)$, then $$\int_{B(x, r)\cap \Omega} h^2dx\asymp h(x+r)^2r^n.$$
\end{proposition}

\noindent\textbf{Proof:} See Gyrya and Saloff-Coste \cite{gyryasaloffcoste1}, Proposition 6.4 and Proposition 6.7.\qed\\

\section{Heat Kernel Estimates of the Laplacian in a Lipschitz Domain}

In this section we review the results of Gyrya and Saloff-Coste in \cite{gyryasaloffcoste1} on domains above Lipschitz graphs.\\

\begin{proposition}
Let $\Omega$ be the domain above the graph of a Lipschitz function. Then $\Omega$ is uniform, and hence satisfies the uniform boundary Harnack principle.
\end{proposition}

\noindent\textbf{Proof:} Gyrya and Saloff-Coste \cite{gyryasaloffcoste1}, Proposition 6.6.\qed

\begin{lemma}\label{nondec}
Let $\Omega$ be the domain above the graph of a Lipschitz function. Let $h>0$ be a harmonic profile for $\Omega$. Then if $x,y\in \Omega$ with $x'=y'$ and $x_n<y_n$, we have $h(x)\leq h(y)$.
\end{lemma}

\noindent\textbf{Proof:} Gyrya and Saloff-Coste \cite{gyryasaloffcoste1}, Proposition 6.4.\qed

\begin{theorem}
Let $\Omega$ be the domain above the graph of a Lipschitz function. Let $h>0$ be a harmonic profile for $\Omega$ and let $p_\Omega(t,x,y)$ be the heat kernel of the Laplacian in $\Omega$ with Dirichlet boundary conditions. Then there exist constants $c_1, c_2, c_3, c_4>0$ such that for all $t>0$ and all $x,y\in \Omega$ we have $$\frac{c_1h(x)h(y)}{t^{n/2}h(x+\sqrt{t})h(y+\sqrt{t})}e^{-c_2\frac{\abs{x-y}^2}{t}}\leq p_\Omega(t,x,y)\leq\frac{c_3h(x)h(y)}{t^{n/2}h(x+\sqrt{t})h(y+\sqrt{t})}e^{-c_4\frac{\abs{x-y}^2}{t}}.$$
\end{theorem}

\noindent\textbf{Proof:} Gyrya and Saloff-Coste \cite{gyryasaloffcoste1}, Corollary 6.14.\qed

\begin{corollary}\label{greenhke}
Let $\Omega$ be the domain above the graph of a Lipschitz function. Let $h>0$ be a harmonic profile for $\Omega$ and let $G_\Omega(x,y)$ be the Green's function of the Laplacian in $\Omega$. Then there exist constants $c, C>0$ such that for all $x,y\in \Omega$, $$c\int_{\abs{x-y}^2}^\infty \frac{h(x)h(y)}{t^{n/2}h(y+\sqrt{t})^2}dt\leq G_\Omega(x,y)\leq C\int_{\abs{x-y}^2}^\infty \frac{h(x)h(y)}{t^{n/2}h(y+\sqrt{t})^2}dt.$$
\end{corollary}

\noindent\textbf{Proof:} Gyrya and Saloff-Coste \cite{gyryasaloffcoste1}, Theorem 5.13, combined with the fact that the instrinsic distance in $\Omega$ is comparable to Euclidean distance.\qed\\

Corollary \ref{greenhke} has the following consequence, known as the $3G$ Principle.

\begin{proposition}\label{3g}
Let $\Omega$ be the domain above the graph of a Lipschitz function. Let $h>0$ be a harmonic profile for $\Omega$ and let $G_\Omega(x,y)$ be the Green's function of the Laplacian in $\Omega$. Then there exists a constant $C>0$ such that for all $x,y,z\in \Omega$ distinct, we have $$\frac{G_\Omega(x,z)G_\Omega(z,y)}{G_\Omega(x,y)}\leq C\left(\frac{h(z)}{h(x)}G_\Omega(x,z)+\frac{h(z)}{h(y)}G_\Omega(z,y)\right).$$
\end{proposition}

\noindent\textbf{Proof:} For $x,y\in \Omega$ distinct let $G_h(x,y)=\frac{G_\Omega(x,y)}{h(x)h(y)}$. By cross-multiplying it suffices to show that for $x,y,z\in \Omega$ distinct, $$\frac{G_h(x,z)G_h(z,y)}{G_h(x,y)}\leq C\left(G_h(x,z)+G_h(z,y)\right).$$ Thus let $x,y,z\in \Omega$ distinct. By the triangle inequality, $\abs{x-y}\leq \abs{x-z}+\abs{z-y}$ and therefore either $\abs{x-y}\leq 2\abs{x-z}$ or $\abs{x-y}\leq 2\abs{z-y}$. Assume that $\abs{x-y}\leq 2\abs{x-z}$. Using Lemma , we have that $$G_h(x,y)\asymp \int_{\abs{x-y}^2}^\infty \frac{dt}{t^{n/2}h(x+\sqrt{t})^2}\geq \int_{4\abs{x-z}^2}^\infty \frac{dt}{t^{n/2}h(x+\sqrt{t})^2}\asymp G_h(x,z),$$ the last relation obtained by the using change of variables $s=\sqrt{2}t$ and applying volume doubling. Therefore in this case we have $G_h(x,z)\preceq G_h(x,y)$ and so $$\frac{G_h(x,z)G_h(z,y)}{G_h(x,y)}\preceq G_h(z,y).$$ If instead $\abs{x-y}\leq 2\abs{z-y}$, the proof is the same except to first use that $G_h(x,y)=G_h(y,x)$. \qed\\

The following lemma and its corollary concern domains in dimension $3$ and higher. We will not need them in this paper, but they may be of independent interest.

\begin{lemma}\label{classical3glemma}
Let $\Omega\subseteq \Real^N$ be the domain above the graph of a Lipschitz function with $N\geq 3$. Let $h>0$ be a harmonic profile for $\Omega$ and let $G_\Omega(x,y)$ be the Green's function of the Laplacian in $\Omega$. Then there exists $C>0$ such that for all $x,y\in \Omega$ distinct, we have $$\frac{h(y)}{h(x)}G_\Omega(x,y)\leq \frac{C}{\abs{x-y}^{N-2}}.$$
\end{lemma}

\noindent\textbf{Proof:} Let $h>0$ be a harmonic profile for $\Omega$. Using Corollary \ref{greenhke} as well as Lemma \ref{nondec}, we have that $$\frac{h(y)}{h(x)}G_\Omega(x,y)\preceq \int_{\abs{x-y}^2}^\infty \frac{h(y)^2dt}{h(y+\sqrt{t})^2t^{N/2}}\leq \int_{\abs{x-y}^2}^\infty \frac{dt}{t^{N/2}}\asymp \frac{1}{\abs{x-y}^{N-2}}.$$\qed

\begin{corollary}
(Classical $3G$ Principle) Let $\Omega\subseteq \Real^N$ be the domain above the graph of a Lipschitz function with $N\geq 3$. Then there exists $C>0$ such that for all $x,y,z\in \Omega$ distinct, $$\frac{G_\Omega(x,z)G_\Omega(z,y)}{G_\Omega(x,y)}\leq C\left(\frac{1}{\abs{x-z}^{N-2}}+\frac{1}{\abs{z-y}^{N-2}}\right).$$
\end{corollary}

\noindent\textbf{Proof:} Combine Lemma \ref{classical3glemma} with Proposition \ref{3g}.\qed

\section{Estimates in a Lipschitz Domain}

We begin this section with some elementary observations in a Lipschitz domain. Given a domain $\Omega$ and $x\in \Omega$, and let $\delta(x)=d(x, \partial \Omega)$ be the distance from $x$ to the boundary of $\Omega$. Also let $\brac{x}=1+\delta(x)$.

\begin{lemma}\label{deltaheightcomp} Let $\fct{f}{\Real^{N-1}}{\Real}$ be Lipschitz and let $L>0$ be such that $\abs{f(x')-f(y')}\leq L\abs{x'-y'}$ for all $x', y'\in \Real^{N-1}$. Let $\Omega=\set{(x', x_N)\in \Real^N\colon x_N>f(x')}$. Then for all $x=(x', x_N)\in \Omega$,
$$\delta(x)\leq x_N-f(x') \leq (2+L)\delta(x).$$
\end{lemma}

\noindent\textbf{Proof:} Let $x=(x_1, \dots, x_N)\in \Omega$. Since $(x_1, \dots, x_{N-1}, f(x_1, \dots, x_{N-1}))\in \partial \Omega$ we have that $\delta(x)\leq x_N-f(x)$.\\

Next, let $y\in \partial \Omega$ be such that $d(x,y)=\abs{x-y}=\delta(x)$. Since $y$ must lie on the graph of $f$ we can write $y=(y_1, \dots, y_{N-1}, f(y_1, \dots, y_{N-1}))$. Since the distance from $x$ to $(x_1, \dots, x_{N-1}, f(x))$ is at least $\delta(x)$ while $\abs{x'-y'}\leq \abs{x-y}=\delta(x)$, we have that $\abs{f(x')-f(y')}\leq L\abs{x'-y'}\leq L\delta(x)$. Thus $$x_N-f(x)=\abs{x-(x', f(x'))}\leq \abs{x-y}+\abs{y-(x', f(x'))}$$ $$\leq \delta(x)+\abs{x'-y'}+\abs{f(x')-f(y')}\leq (2+L)\delta(x).$$\qed

\begin{lemma}\label{deltaheightcomp2}
Let $\fct{f}{\Real^{N-1}}{\Real}$ be Lipschitz with Lipschitz constant $L>0$ and assume that $\abs{f}\leq M$. Let $\Omega$ be the domain above the graph of $f$. Then for $x\in \Omega$ with $x_N>2M$ we have $$\frac{\delta(x)}{2}\leq x_N\leq 2(2+L)\delta(x).$$
\end{lemma}

\noindent\textbf{Proof:} Using the previous lemma, we have that $\delta(x)\leq x_N-f(x')\leq x_N+M\leq 2x_N$.\\

Also, $x_N\leq 2(x_n-M)\leq 2(x_N-f(x'))\leq 2(2+L)\delta(x)$.\qed\\

Let $M$ be such that $\abs{f}\leq M$. Therefore if we define $\Omega_-=\set{(y_1, \dots, y_N)\colon y_N>M}$ and $\Omega_+=\set{(y_1, \dots, y_N)\colon y_N>-M}$ we see that $\Omega_+\subseteq \Omega\subseteq \Omega_-$. Let $h_+$ and $h_-$ be the harmonic profiles of $\Omega_+$ and $\Omega_-$ respectively. They take a very simple form: $h_+(x)=x_N+M$ and $h_-(x)=x_N-M$.\\

\begin{lemma}\label{hpluscomparison}
Let $\Omega$ be the domain above the graph of a bounded Lipschitz function. Let $h_+$, $h_-$, $\Omega_+$, $\Omega_-$ be as defined above. Let $h>0$ be a harmonic profile on $\Omega$. On $h\preceq h_+$ on $\Omega$ and $h_-\preceq h$ on $\Omega_-$.
\end{lemma}

\noindent\textbf{Proof:} We start by proving that $h\leq h_+$ on $\Omega$. Note that $h\leq h_+$ on $\partial \Omega$.\\ 

Let $g_+$ be the solution in $\Omega$ to the Dirichlet problem with boundary data $g_+=h_+$ on $\partial \Omega$. Then by the general theory of the Dirichlet problem, e.g. Classical Potential Theory \cite{armitagegardiner1}, $g_+\leq h_+$ on $\Omega$ and since $h_+$ is bounded on $\partial \Omega$, $g_+$ is bounded on all of $\Omega$. So $h_+-g_+$ is a positive harmonic function on $\Omega$ that vanishes at $\partial \Omega$. By uniqueness of the profile, there exists $c>0$ such that $h_+-g_+=ch$. Hence $h=\frac{1}{c}\left(h_+-g_+\right)\preceq h_+$.\\

The argument that $h_-\preceq h$ on $\Omega_-$ is similar, using that $h$ is bounded on $\partial \Omega_-$ by the previous part.\qed\\

\begin{proposition}\label{3gbound}
Let $\Omega$ be the domain above the graph of a bounded Lipschitz function. Let $W\in C^\infty(\Omega)$ be such that there exist $c, \epsilon>0$ with $$\abs{W(x)}\leq c\brac{x}^{-(2+\epsilon)}$$ for all $x\in \Omega$, where $\brac{x}=1+\delta(x)$. Then $$\sup_{x\in \Omega}\int_\Omega\frac{h(y)G_\Omega(x,y)}{h(x)}\abs{W(y)}dy<+\infty.$$
\end{proposition}

\noindent\textbf{Proof:} Let $y\in \Omega$ and $W$ be as in the hypothesis. In this proof we extend $W$ to all of $\Real^N$ by defining $W=0$ outside of $\Omega$. Without loss of generality we may assume that $\epsilon<1$. Let $M>0$ be such that $\abs{f}\leq M$. Let $\Omega_+$, $\Omega_-$, $h_+$, and $h_-$ be as in Lemma \ref{hpluscomparison}. Let $C>0$ be such that $h\leq Ch_+$ on $\Omega$ and $h_-\leq Ch$ on $\Omega_-$. Let $y\in \Omega$ and $t>0$. If $t>9M^2$ then $\frac{1}{3}\sqrt{t}\leq y_N+\sqrt{t}-M=h_-(y+\sqrt{t})$.\\

If $y_N>M$ then $h_+(y)=y_N+M\leq 2y_n$. In particular, if $t>9M^2$ and $y_n>M$ then \begin{equation}\label{nicehalfspacebound}\frac{h(y)^2}{h(y+\sqrt{t})^2}\leq \frac{36C^2y_N^2}{t}.\end{equation}

Next, we note that the Green's function $G_\Omega(x,y)$ has the following estimate: \begin{equation}G_\Omega(x,y)\asymp h(x)h(y)\int_{\abs{x-y}^2}^\infty \frac{dt}{h(y+{\sqrt{t}})^2t^{N/2}}.\end{equation} Then $$\int_\Omega \frac{h(y) G_\Omega(x,y)}{h(x)}\abs{W(y)}dy\asymp \int_\Omega \int_{\abs{x-y}^2}^\infty \frac{h(y)^2}{h(y+{\sqrt{t}})^2t^{N/2}}\abs{W(y)}dtdy$$ $$=\int_0^\infty \int_{B(x, \sqrt{t})\cap \Omega} \frac{h(y)^2}{h(y+{\sqrt{t}})^2t^{N/2}}\abs{W(y)}dydt.$$\\ 

 We have $h(y+\sqrt{t})\geq h(y)$ for all $y\in \Omega$ and all $t>0$ by Gyrya and Saloff-Coste \cite{gyryasaloffcoste1}. Thus $$\int_0^{9M^2}\int_{B(x, \sqrt{t})\cap \Omega}\frac{h(y)^2}{h(y+\sqrt{t})^2t^{N/2}}\abs{W(y)}dydt\leq \int_0^{9M^2}\int_{B(x, \sqrt{t})\cap \Omega} \frac{\abs{W(y)}}{t^{N/2}}dydt$$ $$\leq c\int_0^{9M^2}\int_{B(x, \sqrt{t})}\frac{dy}{t^{N/2}}dt\preceq c\cdot 9M^2.$$ Similarly, $$\int_0^{\delta(x)^2/4}\int_{B(x, \sqrt{t})\cap \Omega} \frac{h(y)^2}{h(y+\sqrt{t})^2t^{\frac{N}{2}}}\abs{W(y)}dydt\leq \int_0^{\delta(x)^2/4}\int_{B(x, \sqrt{t})\cap \Omega} \frac{c\brac{x}^{-(2+\epsilon)}}{t^{\frac{N}{2}}}dydt$$ $$\preceq \int_0^{\delta(x)^2/4}\brac{x}^{-(2+\epsilon)}dt\preceq \brac{x}^{-\epsilon}.$$ From these computations we conclude that \begin{equation}\label{firstpartintegral}\sup_{x\in \Omega}\int_0^{\max(9M^2, \delta(x)^2/4)}\int_{B(x, \sqrt{t})\cap \Omega} \frac{h(y)^2}{h(y+\sqrt{t})^2} \abs{W(y)}dydt<+\infty.\end{equation}
 
Next assume that $t>\max(9M^2, \delta(x)^2/4)$. Using $t>9M^2$ we have that

$$\int_{B(x, \sqrt{t})\cap \Omega} \frac{h(y)^2}{h(y+\sqrt{t})^2t^{n/2}}\abs{W(y)}dy\leq 9C^2\int_{B(x, \sqrt{t})\cap \Omega} \frac{(y_N+M)^2}{t^{1+\frac{N}{2}}}\abs{W(y)}dy$$ $$\leq \frac{9C^2}{t^{1+\frac{N}{2}}}\int_{B(x, \sqrt{t})\cap \Omega} (y_N+M)^2\abs{W(y)}dy.$$

 Note that for $y\in \Omega$ with $y_N>2M$, by Lemma \ref{deltaheightcomp} we have $y_N\leq 2(2+L)\delta(y)$ and thus $\abs{W(y)}\preceq y_N^{-(2+\epsilon)}$. Hence letting $Q(x, \sqrt{t})=B(x', \sqrt{t})\times (x_N-\sqrt{t}, x_N+\sqrt{t})$ we note that $B(x, \sqrt{t})\subseteq Q(x, \sqrt{t})$ and thus $$\int_{B(x, \sqrt{t})\cap \Omega} (y_N+M)^2\abs{W(y)}dy \leq \int_{Q(x, \sqrt{t})\cap \Omega} (y_N+M)^2\abs{W(y)}dy$$ $$\leq\int_{B(x', \sqrt{t})}\int_{x_N-\sqrt{t}}^{x_N+\sqrt{t}}(y_N+M)^2\abs{W(y',y_N)}dy_Ndy'.$$ Next, $$\int_{x_N-\sqrt{t}}^{x_N+\sqrt{t}} (y_N+M)^2\abs{W(y', y_N)}dy_N$$ $$=\int_{x_N-\sqrt{t}}^{\max(x_N-\sqrt{t}, 2M)}(y_N+M)^2\abs{W(y', y_N)}dy_N+\int_{\max(x_N-\sqrt{t}, 2M)}^{x_N+\sqrt{t}}(y_N+M)^2\abs{W(y', y_N)}dy_N$$ $$\preceq \int_{-2M}^{2M} 9M^2\abs{W(y', y_N)}dy_N+\int_{\max(x_N-\sqrt{t}, 2M)}^{x_N+\sqrt{t}} (y_N)^2(y_N)^{-(2+\epsilon)}dy_n$$ $$\leq 36M^3 \left(\sup_{\Omega_-\cap \Omega_+} \abs{W}\right)+\int_{\max(x_N-\sqrt{t}, 2M)}^{x_N+\sqrt{t}}y_N^{-\epsilon}dy_N$$ $$\leq 36M^3\left(\sup_{\Omega_-\cap \Omega_+} \abs{W}\right)+\frac{(x_N+\sqrt{t})^{1-\epsilon}}{1-\epsilon}.$$

Here we utilize Lemma \ref{deltaheightcomp2} and $t>\delta(x)^2/4$ to see that $x_N+\sqrt{t}\leq 2(2+L)\delta(x)+\sqrt{t}<(3+L)\sqrt{t}$. Hence \begin{equation}36M^3\left(\sup_{\Omega_-\cap \Omega_+} \abs{W}\right)+\frac{(x_N+\sqrt{t})^{1-\epsilon}}{1-\epsilon}\preceq t^{\frac{1-\epsilon}{2}}.\end{equation}

We therefore get that $$\int_{B(x, \sqrt{t})\cap \Omega}(y_N+M)^2\abs{W(y)}dy\leq \int_{B(x', \sqrt{t})}\int_{x_N-\sqrt{t}}^{x_N+\sqrt{t}}(y_N+M)^2\abs{W(y', y_N)}dy_Ndy'$$ $$\preceq \int_{B(x', \sqrt{t})}t^{\frac{1-\epsilon}{2}}dy\asymp t^{\frac{N-1}{2}}\cdot t^{\frac{1-\epsilon}{2}}=t^{\frac{N-\epsilon}{2}}.$$

Thus we have $$\int_{\max(9M^2, \delta(x)^2/4)}^\infty \int_{B(x, \sqrt{t})\cap \Omega} \frac{h(y)^2}{h(y+\sqrt{t})^2t^{\frac{N}{2}}}\abs{W(y)}dydt\preceq \int_{\max(9M^2, \delta(x)^2/4)}^\infty \frac{t^{\frac{N-\epsilon}{2}}}{t^{1+\frac{N}{2}}}dt$$ \begin{equation}\label{tailbound}=\int_{\max(9M^2, \delta(x)^2/4)}^\infty \frac{dt}{t^{1+\frac{\epsilon}{2}}}\preceq (\max(9M^2, \delta(x)^2/4))^{-\frac{\epsilon}{2}}.\end{equation} We thus conclude that \begin{equation}\label{secondpartintegral} \sup_{x\in \Omega}\int_{\max(9M^2, \delta(x)^2/4)}^\infty \int_{B(x, \sqrt{t})\cap \Omega} \frac{h(y)^2}{h(y+\sqrt{t})^2t^{N/2}}\abs{W(y)}dy<+\infty.
\end{equation}

Combining \ref{firstpartintegral} with \ref{secondpartintegral} we see that $$\sup_{x\in \Omega} \int_\Omega \frac{h(y)G_\Omega(x,y)}{h(x)}\abs{W(y)}dy\asymp \sup_{x\in\Omega}\int_0^\infty \int_{B(x, \sqrt{t})\cap \Omega} \frac{h(y)^2}{h(y+\sqrt{t})^2t^{N/2}}\abs{W(y)}dydt$$ $$\leq \sup_{x\in \Omega} \int_0^{\max(9M^2, \delta(x)^2/4)} \int_{B(x, \sqrt{t})\cap \Omega} \frac{h(y)^2}{h(y+\sqrt{t})^2t^{N/2}}\abs{W(y)}dydt$$ $$+\sup_{x\in \Omega} \int_{\max(9M^2, \delta(x)^2/4)}^\infty \int_{B(x, \sqrt{t})\cap \Omega} \frac{h(y)^2}{h(y+\sqrt{t})^2t^{N/2}}\abs{W(y)}dydt$$ $$<+\infty.$$\qed

\begin{corollary}\label{3gbound2}
Let $\Omega$ be the domain above the graph of a bounded Lipschitz function. Let $W\in C^\infty(\Omega)$ be such that there exist $c, \epsilon>0$ with $$\abs{W(x)}\leq c\brac{x}^{-(2+\epsilon)}$$ for all $x\in \Omega$, where $\brac{x}=1+\delta(x)$. Then $$\sup_{x,y\in \Omega}\int_\Omega \frac{G_\Omega(x,z)G_\Omega(z,y)\abs{W(z)}}{G_\Omega(x,y)}dz<+\infty.$$
\end{corollary}

\noindent\textbf{Proof:} From Proposition \ref{3g}, let $h>0$ be a harmonic profile for $\Omega$ and let $C>0$ be such that for all $x,y,z\in \Omega$ distinct, we have $$\frac{G_\Omega(x,z)G_\Omega(z,y)}{G_\Omega(x,y)}\leq C\left(\frac{h(z)}{h(x)}G_\Omega(x,z)+\frac{h(z)}{h(y)}G_\Omega(z,y)\right).$$ Then using Proposition \ref{3gbound}, we get that $$\sup_{x,y\in \Omega} \int_\Omega \frac{G_\Omega(x,z)G_\Omega(z,y)\abs{W(z)}}{G_\Omega(x,y)}dz$$ $$\leq C\left(\sup_{x\in \Omega}\int_\Omega\frac{h(z)G_\Omega(x,z)}{h(x)}\abs{W(z)}dz+\sup_{y\in \Omega} \int_\Omega \frac{h(z)G_\Omega(z,y)}{h(y)}\abs{W(z)}dz\right)$$ $$<+\infty.$$ \qed

\section{Conditional Gaugeability}

In this section we discuss conditional Brownian motion in $\Omega$. For a textbook presentation of these ideas, see Chung and Zhao \cite{chungzhao1} Chapter 5.\\

Let $(X_t)_{t>0}$ denote Brownian motion in $\Omega$ with killing at the boundary. The Markov transition kernel of $(X_t)_{t>0}$ is $p^\Omega(t,x,z)$, the heat kernel of the Laplacian in $\Omega$ with Dirichlet boundary conditions. Fix $y\in \Omega$. We now consider a new Markov transition kernel given by $$p_y^\Omega(t,x,z)=\frac{G_\Omega(z,y)p^\Omega(t,x,z)}{G_\Omega(x,y)},$$ where $x,z\neq y$. The Markov process in $\Omega\setminus\set{y}$ defined by $p_y^\Omega(t,x,z)$ is Brownian motion conditioned to end at $y$, denoted $(X^{\cdot, y}_t)_{t>0}$. Let $\zeta^y$ denote the lifetime of this process. Let $E^x_y$ denote expectation with respect to this process started at $x\in \Omega$.\\

The following lemma can be found e.g. in the proof of Lemma 3 in Zhao \cite{zhao1}. We reproduce the argument here.

\begin{lemma}\label{3gprobbound}
Let $\Omega$ be the domain above the graph of a bounded Lipschitz function and let $W\in C^\infty(M)$ be bounded. Let $x,y\in \Omega$ be distinct. Then $$E^x_y\left[\int_0^{\zeta^y} \abs{W(X^{\cdot, y}_s)}ds\right]\leq\int_\Omega \frac{G_\Omega(x,z)G_\Omega(z,y)}{G_\Omega(x,y)}\abs{W(z)}dz.$$
\end{lemma}

\noindent\textbf{Proof:} For $x,y\in \Omega$ distinct, by the Fubini-Tonelli Theorem we have that $$E^x_y\left[\int_0^{\zeta^y}\abs{W(X^{\cdot, y}_s)}ds\right]=\int_0^\infty E^x_y\left[s<\zeta^y; \abs{W(X^{\cdot, y}_s)}\right]ds$$ $$=\frac{1}{G_\Omega(x,y)}\int_0^\infty E^x\left[s<\zeta^y; G_\Omega(X_s, y)\abs{W(X_s)}\right]ds$$ $$\leq\frac{1}{G_\Omega(x,y)}\int_0^\infty \int_\Omega p^\Omega(s,z,y)G_\Omega(z,y)\abs{W(z)}dzdt$$ $$=\int_\Omega\frac{G_\Omega(x,z)G_\Omega(z,y)}{G_\Omega(x,y)}\abs{W(z)}dz.$$\qed\\

Suppose in particular that $W\geq 0$ and $a=\sup_{x,y\in \Omega}\int_\Omega \frac{G_\Omega(x,z)G_\Omega(z,y)}{G_\Omega(x,y)}W(z)dz<+\infty$. Then by Lemma \ref{3gprobbound} and Jensen's inequality we have that $$1\geq E^x_y\left[\exp\left(-\int_0^{\zeta^y} W(X^{\cdot, y}_s)ds\right)\right]\geq \exp\left(-E^x_y\left[\int_0^{\zeta^y} W(X^{\cdot, y}_s)ds\right]\right)\geq e^{-a}.$$

\begin{definition}\label{conditionalgaugedef}
Let $\fct{W}{\Omega}{\Real}$. For $t>0$ we define $e_W(t)=-\int_0^t W(X^{\cdot, y}_s)ds$. We define the \textbf{conditional gauge} with respect to $W$ as the function $$E^x_y[e_W(\zeta^y)].$$ We say that $W$ is \textbf{conditionally gaugeable} if $$\sup_{x,y\in \Omega}E^x_y[e_W(\zeta^y)]<+\infty.$$
\end{definition}

\begin{proposition}\label{greenratio}
Assume that $W$ is conditionally gaugeable. Then $\Delta+W$ has a Green's function, denoted $G_\Omega^W(x,y)$, and in particular, for $x,y\in \Omega$ distinct, $$G_\Omega^W(x,y)=G_\Omega(x,y)E^x_y[e_W(\zeta^y)].$$
\end{proposition}

\noindent\textbf{Proof:} See in particular Chen \cite{chen1}, Lemma 3.5, Theorem 3.6, and the comment after the proof of Theorem 3.6.\qed\\

\begin{corollary}\label{greenwbound}
Let $\Omega$ be the domain above the graph of a bounded Lipschitz function. Let $W\in C^\infty(\Omega)$ be such that there exist $c, \epsilon>0$ with $$0\leq W(x)\leq c\brac{x}^{-(2+\epsilon)}$$ for all $x\in \Omega$, where $\brac{x}=1+\delta(x)=1+d(x, \partial \Omega)$. Then the Green's function $G_\Omega^W(x,y)$ for $\Delta+W$ in $\Omega$ exists and there exists $C>0$ such that for all $x,y\in \Omega$ distinct we have $$\frac{1}{C}G_\Omega(x,y)\leq G^W_\Omega(x,y)\leq G_\Omega(x,y).$$ 
\end{corollary}

\noindent\textbf{Proof:} Follows from Proposition \ref{greenratio} plus the Jensen's inequality remark above Definition \ref{conditionalgaugedef}.\qed\\

The main useful property of $G_\Omega^W(x,y)$, besides that given in Proposition \ref{greenratio}, is that for $x\in \Omega$ fixed, the function $u=G_\Omega^W(x,\cdot)$ satisfies $(\Delta+W)u=0$ in $\Omega\setminus\set{x}$. We now introduce the classical Harnack principle for positive solutions to $(\Delta+W)u=0$.  

\begin{lemma}\label{schrodingerharnack}
Let $D$ be a domain in $\Real^N$ and let $W$ be a bounded Borel function on $D$. Let $(u_n)_{n\in \N}$ be a sequence a functions on $D$ such that for each $n\in \N$, $u_n>0$ on $D$ and $(\Delta+W)u_n=0$ in $D$. Assume also that the sequence $(u_n)_{n\in \N}$ is bounded uniformly on each compact subset of $D$. Then there exists a function $v$ on $D$ such that $(\Delta+W)v=0$ and a subsequence $(u_{n_k})$ such that $u_{n_k}\rightarrow v$ locally uniformly on $D$.
\end{lemma}

\noindent\textbf{Proof:} The sequence $(u_n)_{n\in \N}$ is equicontinuous by e.g. Theorem 8.22 in Gilbarg and Trudinger's text \cite{gilbargtrudinger1}. The conclusion follows from the Arzela-Ascoli Theorem.\qed\\

\begin{theorem}\label{profilebound}
Let $\Omega$ be the domain above the graph of a bounded Lipschitz function. Let $W\in C^\infty(\Omega)$ be such that there exist $c, \epsilon>0$ with $$\abs{W(x)}\leq c\brac{x}^{-(2+\epsilon)}$$ for all $x\in \Omega$, where $\brac{x}=1+\delta(x)=1+d(x, \partial \Omega)$. Then a profile $h^W>0$ for $\Delta+W$ in $\Omega$ exists and, if $h>0$ is a harmonic profile in $\Omega$, there exists a constant $C>0$ such that for all $x\in \Omega$, $$\frac{1}{C}h(x)\leq h^W(x)\leq Ch(x).$$
\end{theorem}

\noindent\textbf{Proof:} We follow the construction of the harmonic profile in $\Omega$ given in \cite{gyryasaloffcoste1} using Green's function, now with the tools to compare the classical Green's function to that of $\Delta+W$. Let constants $A_0, A_1>1$ be as in the statement of the Boundary Harnack Principle. Let $C>0$ be as in Corollary \ref{greenwbound}.\\

Fix $o\in \Omega$ and fix $\xi\in \partial \Omega$. Let $(r_n)_{n\in \N}$ be an increasing sequence of positive radii with $r_n\rightarrow +\infty$. Let $(x_n)_{n\in \N}$ be a sequence in $\Omega$ such that for all $n\in \N$, $x_n\notin B(\xi, A_0r_n)$. Fix a distinguished point $o\in \Omega$ distinct from each $x_n$. Define $$h_n(x):=\frac{G_\Omega(x, x_n)}{G_\Omega(o, x_n)},$$ $$h_n^W(x):=\frac{G_\Omega^W(x, x_n)}{G_\Omega^W(o, x_n)}.$$ Note that $h_n(o)=h_n^W(o)=1$ for all $n$ and by Corollary \ref{greenwbound}, we have \begin{equation}\frac{1}{C}h_n(x)\leq h_n^W(x)\leq Ch_n(x)\end{equation} for all $x\in \Omega\setminus \set{x_n}$ and all $n$. Fix $n_0\in \N$ and let $m,n \in \N$ with $m,n\geq n_0$. Since $x_n\notin B(\xi, A_0r_{n_0})$, the function $G_\Omega(\cdot, x_n)$ is harmonic in $B(\xi, A_0r_{n_0})\cap \Omega$ with Dirichlet boundary conditions along $\partial \Omega$. Similarly for $m$. Thus by the Boundary Harnack Principle \ref{boundaryharnack}, we have that for any $x\in B(\xi, r_{n_0})\cap \Omega$, $$h_n(x)=\frac{G_\Omega(x, x_n)}{G_\Omega(o, x_n)}\leq A_1\frac{G_\Omega(x, x_m)}{G_\Omega(o, x_m)}=A_1h_m(x).$$ By symmetry we also have $h_m(x)\leq A_1 h_n(x)$. We conclude that the sequence $(h_n)_{n\geq n_0}$ is uniformly bounded on each compact subset of $B(\xi, r_{n_0})\cap \Omega$, and hence the full sequence $(h_n)_{n\in \N}$ is, for each compact subset $K\subseteq \Omega$, eventually uniformly bounded on $K$.\\

 Therefore by the classical Harnack principle for harmonic functions as well as passing to a subsequence, we may assume that the sequence $(h_n)_{n\in \N}$ converges uniformly on compact sets. By $\frac{1}{C}h_n\leq h_n^W\leq Ch_n$, we know that $(h_n^W)_{n\in \N}$ is uniformly bounded on each compact subset of $\Omega$. So by Lemma \ref{schrodingerharnack}, we may pass to a subsequence again and assume that $(h_n^W)_{n\in \N}$ converges uniformly on compact subsets of $\Omega$. Let $$h(x)=\lim_{n\rightarrow\infty} h_n(x),$$ $$h^W(x)=\lim_{n\rightarrow \infty} h_n^W(x).$$ It is routine to verify that $h$ is harmonic and $(\Delta+W)h^W=0$. Clearly $h>0$. Also, since $\frac{1}{C}h_n\leq h_n^W\leq Ch_n$ for all $n$, we must have $$\frac{1}{C}h\leq h^W\leq Ch.$$ By Theorem 4.16 in Gyrya and Saloff-Coste \cite{gyryasaloffcoste1}, $h$ satisfies Dirichlet boundary conditions along $\partial \Omega$, i.e. $h$ is a harmonic profile for $\Omega$. By $\frac{1}{C}h\leq h^W\leq Ch$, $h^W$ also satisfies Dirichlet boundary conditions and thus is a profile for $\Delta+W$ in $\Omega$.\qed

\section{Dirichlet Forms and Doob's $h$-Transform}

Equipped with Theorem \ref{profilebound}, we are ready to study the heat kernel of $\Delta+W$ in $\Omega$. First we consider the closed form on $L^2(\Omega, dx)$ given by $$\mc{E}^W(f,f)=\int_\Omega \left(\abs{\grad f}^2 + Wf^2\right)dx$$ with domain $\mc{F}=W^{1,2}_0(\Omega, dx)$. The fact that we used $W^{1,2}_0(\Omega, dx)$ as the domain instead of $W^{1,2}(\Omega, dx)$ corresponds to Dirichlet boundary conditions.\\

Let $h^W>0$ be a profile for $\Delta+W$ as in Theorem \ref{profilebound}, and let $d\mu=(h^W)^2dx$. Consider the strictly local Dirichlet form on $L^2(\Omega, d\mu)$ given by $$\mc{E}_\mu(f,f)=\int_\Omega \abs{\grad f}^2 d\mu.$$ Together with the domain $$\mc{F}^\mu=\set{f\in L^2(\Omega, d\mu)\cap W^{1,2}_{\textrm{loc}}(\Omega, dx)\colon \int_\Omega \abs{\grad f}^2d\mu<+\infty}$$ this is a closed form.\\

Next consider the unitary map $\fct{H}{L^2(\Omega, d\mu)}{L^2(\Omega, dx)}$ given by $H(f)=h^Wf$. This yields a transformation of the form $\mc{E}^W$, allowing us to define $$\mc{E}_H(f,f)=\mc{E}^W(hf, hf)$$ with domain $H^{-1}\mc{F}$. This is a strictly local Dirichlet form on $L^2(\Omega, d\mu)$.

\begin{proposition}\label{formscoincide}
The forms $\mc{E}_H$ and $\mc{E}_\mu$ coincide, i.e. $H^{-1}\mc{F}=\mc{F}_\mu$ and $\mc{E}_H=\mc{E}_\mu$ on their shared domain.
\end{proposition}

\noindent \textbf{Proof:} See Gyrya and Saloff-Coste \cite{gyryasaloffcoste1}, Proposition 5.7 and Theorem 5.9. In that presentation $W=0$, i.e. the profile is harmonic, but this does not affect any of the proofs.\qed\\

Since the forms $\mc{E}_H$ and $\mc{E}_\mu$ coincide, so do their infinitesimal generators, and hence so do their heat kernels. We let $p_\mu(t,x,y)$ denote the shared heat kernel. We have the following useful lemma that relates $p_\mu(t,x,y)$ to the heat kernel of $\Delta+W$, which we refer to as $p^W(t,x,y)$.

\begin{lemma}\label{htransnicelemma}
Let $p_\mu(t,x,y)$ and $p^W(t,x,y)$ be as above. Then for all $t>0$ and all $x,y\in \Omega$, we have $$p^W(t,x,y)=h^W(x)h^W(y)p_\mu(t,x,y).$$
\end{lemma}

\noindent\textbf{Proof:} Follows from Lemma \ref{formscoincide} as well as Gyrya and Saloff-Coste \cite{gyryasaloffcoste1}, Lemma 5.6.\qed\\

\subsection{Heat Kernel Estimates: Main Results}

\begin{theorem}\label{htransheatkernel}
Let $\Delta_\mu$ be the infinitesimal generator of the form $\mc{E}_\mu$. Then the heat kernel of $\Delta_H$ satisfies the following estimate: there exist $c_1,c_2,c_3,c_4>0$ such that for all $t>0$ and all $x,y\in\Omega$ we have $$\frac{c_1}{\mu(B(x,\sqrt{t})\cap \Omega)}e^{-c_2\frac{d(x,y)^2}{t}}\leq p_\mu(t,x,y)\leq \frac{c_3}{\mu(B(x,\sqrt{t})\cap \Omega)}e^{-c_4\frac{d(x,y)^2}{t}}.$$
\end{theorem}

\noindent\textbf{Proof:} This follows from Gyrya and Saloff-Coste \cite{gyryasaloffcoste1}, Theorem 3.34. It is also worth mentioning the series of papers by Sturm \cite{sturm1}, \cite{sturm2}, \cite{sturm3} which the result in the previous sentence uses. \qed\\

The following is our main result.

\begin{theorem}\label{maithm}
Let $\Omega\subseteq \Real^N$ be the domain above the graph of a bounded Lipschitz function and let $W\in C^\infty(\Omega)$ be such that there exist $c, \epsilon>0$ with $$0\leq W(x)\leq c\brac{x}^{-(2+\epsilon)}$$ for all $x\in \Omega$, where $\brac{x}=1+\delta(x)=1+d(x, \partial \Omega)$. Let $h>0$ be a harmonic profile on $\Omega$. If $p^W(t,x,y)$ denotes the Dirichlet heat kernel of $\Delta+W$ on $\Omega$, then there exist constants $c_1, c_2, c_3, c_4>0$ such that for all $t>0$ and $x,y\in \Omega$ we have $$\frac{c_1h(x)h(y)}{h(x+\sqrt{t})h(y+\sqrt{t})t^{N/2}}e^{-c_2\frac{d(x,y)^2}{t}}\leq p^W(t,x,y)\leq \frac{c_3h(x)h(y)}{h(x+\sqrt{t})h(y+\sqrt{t})t^{N/2}}e^{-c_4\frac{d(x,y)^2}{t}}.$$
\end{theorem}

\noindent\textbf{Proof:} Combine Theorem \ref{htransheatkernel} with Theorem \ref{profilebound}, Lemma \ref{htransnicelemma}, and Lemma \ref{volumeest}.\qed\\

\noindent\textbf{Remark:} If the Lipschitz function $f$ defining $\Omega$ satisfies $\abs{f}\leq M$, then for $x,y\in \Omega$ with $x_N, y_N>2M$, we have $h(x)\asymp \delta(x)$. It follows that, using the same notation as Theorem \ref{maithm}, for any $t>0$ and any $x,y\in \Omega$ with $x_N, y_N>2M$, $$\min\left(\frac{\delta(x)\delta(y)}{t}, 1\right)\frac{1}{t^{N/2}}e^{-c_2\frac{d(x,y)^2}{t}}\preceq p^W(t,x,y)\preceq \min\left(\frac{\delta(x)\delta(y)}{t}, 1\right)\frac{1}{t^{N/2}}e^{-c_4\frac{d(x,y)^2}{t}}.$$

Now assume that $f$ is bounded and $C^{1,1}$, in the sense that $f\in C^1$ and $\grad f$ is Lipschitz with uniform Lipschitz constant. Lemma 2.2 in Song \cite{song1} implies that the harmonic profile $h$ for the domain $\Omega$ above the graph of $f$ satisfies $h(x)\asymp \delta(x)$. It follows that there exist $c_1, c_2, c_3, c_4>0$ such that for all $t>0$ and all $x,y\in \Omega$, $$\min\left(\frac{\delta(x)\delta(y)}{t}, 1\right)\frac{c_1}{t^{N/2}}e^{-c_2\frac{d(x,y)^2}{t}}\leq p^W(t,x,y)\leq \min\left(\frac{\delta(x)\delta(y)}{t}, 1\right)\frac{c_3}{t^{N/2}}e^{-c_4\frac{d(x,y)^2}{t}}.$$

\subsection{Obstacles to the Unbounded Case: Example in $\Real^2$}

In this section we show some possible obstacles to a proof of the analogous result when the Lipschitz function $f$ is not necessarily bounded. Consider the function $\fct{f}{\Real}{\Real}$ given by $f(x_1)=-\abs{x_1}$ and let $\Omega=\set{(x_1, x_2)\colon x_2>f(x_1)}$. Viewing $\Real^2$ as $\C$ with polar coordinates $z=re^{i\theta}$, $r\geq 0$ and $-\pi/2\leq \theta<3\pi/2$, we then have \begin{equation}\Omega=\set{re^{i\theta}\colon r>0, -\pi/4<\theta<5\pi/4}.\end{equation}

With this definition, the harmonic profile of $\Omega$ is then given in polar coordinates as \begin{equation}h(re^{i\theta})=r^{2/3}\sin\left(\frac{\pi}{6}+\frac{2\theta}{3}\right).\end{equation} Notably, when trying to bound the quantity $\frac{h(y)^2}{h(y+\sqrt{t})^2}$, we do not get a tractable expression as in \ref{nicehalfspacebound}. If hypothetically we achieve an upper bound of the form $t^{-2/3}$, then we have the additional obstacle of $t^{-\left(\frac{2}{3}+\frac{\epsilon}{2}\right)}$ not having an integrable tail for small $\epsilon>0$. Compare this to \ref{tailbound}, an important bound in the proof involving an integrable tail.

\section{Acknowledgements}

I would like to thank my advisor, Laurent Saloff-Coste, for helpful comments and support. This work was supported in part by the National Science Foundation Graduate Research Fellowship grants number DGE-2139899 and DGE-1650441.

\bibliographystyle{plain}

\bibliography{LipschitzBib}

\end{document}